\documentclass{amsart}

\usepackage{amssymb}

\newtheorem{lem}{Lemma}[section]
\newtheorem{theo}[lem]{Theorem}
\newtheorem{stat}[lem]{Proposition}
\newtheorem{cons}[lem]{Corollary}

\newtheorem{corollary}[lem]{Corollary}
\newtheorem{problem}[lem]{Problem}

\theoremstyle{definition}

\newtheorem{defn}[lem]{Definition}
\newtheorem{rem}[lem]{Remark}

\newtheorem{remark}[lem]{Remark}

\def\cl{\operatorname{cl}}

\def\clin{\mathrel{\underset{cl}{\subset}}}
\newcommand{\F}{\mathcal {F}}

\newcommand{\U}{\mathcal U}
\newcommand{\V}{\mathcal V}
\newcommand{\W}{\mathcal W}
\newcommand{\A}{\mathcal A}
\newcommand{\Comp}{\mathcal Comp}
\newcommand{\Id}{\mathrm{Id}}
\newcommand{\Fil}{\mathrm{Fil}}
\newcommand{\w}{\omega}
\newcommand{\supp}{\mathrm{supp}}
\newcommand{\IN}{\mathbb N}
\newcommand{\IZ}{\mathbb Z}

\newcommand{\C}{\mathcal C}
\newcommand{\invG}{\overset{\leftrightarrow}{G}}
\newcommand{\inv[1]}{\overset{\leftrightarrow}{#1}}
\newcommand{\gen[1]}{\langle #1\rangle}

\title{Right-topological semigroup operations on inclusion hyperspaces}
\author{Volodymyr Gavrylkiv}
\email{vgavrylkiv@yahoo.com}

\begin{document}
\begin{abstract}We show that for any discrete semigroup $X$ the semigroup operation can be extended to a right-topological semigroup operation on the space $G(X)$ of inclusion hyperspaces on $X$. We detect some important subsemigroups of $G(X)$, study the minimal ideal, the (topological) center, left cancelable elements of $G(X)$, and describe the structure of the semigroups $G(\IZ_n)$ for small numbers $n$.
\end{abstract}
\subjclass{22A15, 54D35}
\maketitle
\tableofcontents

\section*{Introduction}
 After the topological proof  of
Hindman theorem \cite{Hind} given by Galvin and Glazer (unpublished, see 
\cite[p.102]{HS}, \cite{H2}) topological methods become a standard
tool in the modern combinatorics of numbers, see \cite{HS},
\cite{P}. The crucial point is that the semigroup operation
$\ast$ defined on any discrete space $S$ can be extended to a
right-topological semigroup operation on $\beta S$, the Stone-\v
Cech compactification of $S$. The product of two ultrafilters
$\U,\V\in\beta S$ can be found in two steps: firstly for every
element $a\in S$ of the semigroup we extend the left shift
$L_a:S\to S$, $L_a:x\mapsto a\ast x$, to a continuous map $\beta
L_a:\beta S\to\beta S$. In such a way, for every $a\in S$ we
define the product $a*\V=\beta L_a(\V)$. Then, extending the
function $R_\V:S\to\beta S$, $R_\V:a\mapsto a*\V$, to a continuous
map $\beta R_\V:\beta S\to\beta S$, we define the product
$\U\circ\V=\beta R_\V(\U)$. This product can be also defined
directly: this is an ultrafilter with the base
$\bigcup_{x\in U}x\ast V_x$ where $U\in\U$ and $\{V_x\}_{x\in U}\subset\V$.
Endowed with so-extended operation the Stone-\v
Cech compactification $\beta S$ becomes a compact Hausdorff
right-topological semigroup. Because of the compactness the
semigroup $\beta S$ has idempotents, minimal (left) ideals, etc., whose
existence has many important combinatorial consequences.

The Stone-\v Cech compactification $\beta S$ can be considered as a subset
 of the double power-set $\mathcal P(\mathcal P(S))$.  The power-set
  $\mathcal P(X)$ of any set $X$ (in particular, $X=\mathcal P(S)$) carries a natural compact Hausdorff
  topology inherited from the Cantor cube  $\{0,1\}^{X}$
  after identification of each subset $A\subset X$ with its characteristic
   function. The power-set $\mathcal P(X)$ is a complete distributive
 lattice with respect to the operations of union and intersection.

The smallest complete sublattice of $\mathcal P(\mathcal P(S))$
containing $\beta S$ coincides with the space $G(S)$ of inclusion
hyperspaces, a well-studied object in Categorial Topology. By
definition, a family $\mathcal{A}\subset\mathcal P(S)$ of
non-empty subsets of $S$ is called an {\em inclusion hyperspace}
if together with each set $A\in \mathcal{A}$ the family
$\mathcal{A}$ contains all supersets of $A$ in $S$. In \cite{Gav}
it is shown that $G(S)$ is a compact Hausdorff lattice with
respect to the operations of intersection and union.

Our principal observation is that the algebraic operation of the
semigroups $S$ can be extended not only to $\beta S$ but also to
the complete lattice hull $G(S)$ of $\beta S$ in $\mathcal
P(\mathcal P(S))$. Endowed with so-extended operation, the space of
inclusion hyperspaces $G(S)$ becomes a compact Hausdorff
right-topological semigroup containing $\beta S$ as a closed
subsemigroup. Besides $\beta S$, the semigroup $G(S)$ contains many
other important spaces as closed subsemigroups: the superextension
$\lambda S$ of $S$, the space $N_k(S)$ of $k$-linked inclusion
hyperspaces, the space $\Fil(S)$ of filters on $S$ (which contains an isomorphic copy of the global semigroup $\Gamma(S)$ of $S$), etc.

We shall study some properties of the semigroup operation on
$G(S)$ and its interplay with the lattice structure of $G(S)$. We
expect that studying the algebraic structure  of $G(S)$ will
have some combinatorial consequences that cannot be obtained with help of
ultrafilters, see \cite{BGN} for further development of this subject.

\section{Inclusion hyperspaces}

In this section we recall some basic information about inclusion hyperspaces. More detail information can be found in the paper \cite{Gav}.

\subsection{General definition and reduction to the compact case}
For a topological space $X$ by $\exp(X)$ we denote the space of all non-empty closed
 subspaces of $X$ endowed with the Vietoris topology. By an {\em inclusion hyperspace}
  we mean a closed subfamily $\F\subset\exp(X)$ that is {\em monotone} in the sense that
  together with each set $A\in\F$ the family $\F$ contains all closed subsets $B\subset X$
  that contain $A$. By \cite{Gav}, the closure of each monotone family in $\exp(X)$ is an
   inclusion hyperspace. Consequently, each family $\mathcal B\subset\exp(X)$ generates an
    inclusion hyperspace $$\cl_{\exp(X)}\{A\in\exp(X):\exists B\in\mathcal B \mbox{ with }B\subset A\}$$
    denoted by $\gen[\mathcal B]$.\footnote{In \cite{Gav} the inclusion hyperspace $\langle \mathcal B\rangle$ generated by a base $\mathcal B$ is denoted by $\overline{{\uparrow}\mathcal B}$.} In this case $\mathcal B$ is called a {\em base}
    of $\mathcal F=\gen[\mathcal B]$.
An inclusion hyperspace $\gen[x]$ generated by a
singleton $\{x\}$, $x\in X$, is called {\em principal}. 

If $X$ is
discrete, then each monotone family in $\exp(X)$ is an inclusion
hyperspace, see \cite{Gav}.

Denote by $G(X)$
the space of all inclusion hyperspaces with the topology generated
by the subbase
$$
\begin{aligned}
U^+=\;&\{\mathcal{A}\in G(X):\exists B \in\mathcal{A}\mbox{ with }B\subset U\}\mbox{ and }\\
U^{ - }=\;&\{\mathcal{A} \in G(X):\forall B\in\mathcal A\quad B\cap U\ne \emptyset\},
\end{aligned}
$$ where $U$
is open in $X$.

For a $T_1$-space $X$ the map $\eta X:X\to G(X)$, $\eta
X(x)=\{F\clin X:x\in F\}$, is an embedding (see \cite{Gav}), so we
can identify principal inclusion hyperspaces with elements of the
space $X$.

For a $T_1$-space $X$ the space $G(X)$ is Hausdorff if and only if
the space $X$ is normal, see \cite{Gav}, \cite{M2}. In the latter
case the map $$h: G(X)\to G(\beta X),\quad h(\F)=\cl_{\exp(\beta
X)}\{\cl_{\beta X}{F}\,|\,F\in \F\},$$ is a homeomorphism, so we
can identify the space $G(X)$ with the space $G(\beta X)$ of inclusion hyperspaces
over the Stone-\v Cech compactification $\beta X$ of
the normal space $X$, see \cite{M2}. Thus we reduce the study of
inclusion hyperspaces over normal topological spaces to the
compact case where this construction is well-studied.

For a (discrete) $T_1$-space the space $G(X)$ contains a (discrete and) dense subspace $G^\bullet(X)$ consisting of inclusion hyperspaces with finite support. An inclusion hyperspace $\A\in G(X)$ is defined to have {\em finite support in $X$} if $\A=\gen[\mathcal F]$ for some finite family $\F$ of finite subsets of $X$.

An inclusion hyperspace $\F\in G(X)$ on a non-compact space $X$ is called {\em free} if for each compact subset $K\subset X$ and any element $F\in\F$ there is another element $E\in\F$ such that $E\subset F\setminus K$.
By $G^\circ(X)$ we shall denote the subset of $G(X)$ consisting of free inclusion hyperspaces. By \cite{Gav}, for a normal locally compact space $X$ the subset $G^\circ(X)$ is closed in $G(X)$. In the simplest case of a countable discrete space $X=\IN$ free inclusion hyperspaces (called semifilters) on $X=\IN$ have been introduced and intensively studied in \cite{BZ}.

\subsection{Inclusion hyperspaces in the category of compacta} The construction of the space of inclusion hyperspaces
is functorial and monadic in the category $\Comp$ of compact Hausdorff spaces and their continuous
map, see \cite{TZ}. To complete $G$ to a functor on $\Comp$ observe that each continuous map $f:X\to Y$
 between compact Hausdorff spaces induces a continuous map $Gf:G(X)\rightarrow G(Y)$ defined by
$$Gf(\mathcal{A})=\gen[f(\mathcal A)]=\{B \clin Y: B\supset f(A) \text{ for some } A \in
\mathcal{A}\}$$ for $\mathcal{A} \in G(X)$. The map $Gf$ is
well-defined and continuous, and $G$ is a functor in the category
$Comp$ of compact Hausdorff spaces and their continuous maps, see
\cite{TZ}. By Proposition 2.3.2 \cite{TZ}, this functor is weakly
normal in the sense that it is continuous, monomorphic, epimorphic
and preserves intersections, singletons, the empty set and  weight
of infinite compacta.

Since the functor $G$ preserves monomorphisms, for each closed
subspace $A$ of a compact Hausdorff space $X$ the inclusion map
$i:A\to X$ induces a topological embedding $Gi:G(A)\to G(X)$. So
we can identify $G(A)$ with a subspace of $G(X)$. Now for each
inclusion hyperspace $\A\in G(X)$ we can consider the support of
$\A$
$$\supp\A=\cap\{A\clin X:\A\in G(A)\}$$ and conclude that $\A\in G(\supp\A)$ because $G$ preserves
intersections, see \cite[\S2.4]{TZ}.

Next, we consider the monadic properties of the functor $G$.
We recall that a functor $T:\Comp\to\Comp$ is {\em monadic} if it can be completed to a monad ${\mathbb T}=(T,\eta,\mu)$ where $\eta:\Id\to T$ and  $\mu:T^2\to T$ are natural transformations (called the unit and multiplication) such that $\mu\circ T(\mu_X)=\mu\circ \mu_{TX}:T^3X\to TX$ and $\mu\circ \eta_{TX}=\mu\circ T(\eta_X)=\Id_{TX}$ for each compact Hausdorff space $X$, see \cite{TZ}.

For the functor $G$ the unit $\eta:\Id\to G$ has been defined above while the multiplication $\mu=\{\mu_X:G^2X\to G(X)\}$ is defined by the formula
$$\mu _{X}(\Theta)=\cup \{\cap
\mathcal{M}\;|\;\mathcal{M} \in \Theta\}, \;\Theta \in G^2X.$$
By Proposition 3.2.9 of \cite{TZ}, the triple
$\mathbb{G}=(G,\eta ,\mu)$ is a monad in
$\Comp$.

\subsection{Some important subspaces of $G(X)$}\label{subspaces}
The space $G(X)$ of inclusion hyperspaces  contains many interesting subspaces.  Let $X$ be a topological space and $k\ge2$ be a natural number. An inclusion hyperspace $\A\in G(X)$ is defined to be
\begin{itemize}
\item {\em $k$-linked} if $\cap\mathcal F\ne\emptyset$ for any subfamily $\F\subset\A$ with $|\F|\le k$;
\item {\em centered} if $\cap\mathcal F\ne\emptyset$ for any finite subfamily $\F\subset\A$;
\item {\em a filter} if $A_1\cap A_2\in\A$ for all sets $A_1,A_2\in\A$;
\item {\em an ultrafilter} if $\A=\A'$ for any filter $\A'\in G(X)$ containing $\A$;
\item {\em maximal $k$-linked} if $\A=\A'$ for any $k$-linked inclusion hyperspace $\A'\in G(X)$ containing $\A$.
\end{itemize}

By $N_k(X)$, $N_{<\w}(X)$, and $\Fil(X)$ we denote the subsets of $G(X)$ consisting of $k$-linked, centered, and filter inclusion hyperspaces, respectively. Also by $\beta(X)$ and $\lambda_k(X)$ we denote the subsets of $G(X)$ consisting of ultrafilters and maximal k-linked inclusion hyperspaces, respectively.  
The space $\lambda(X)=\lambda_2(X)$ is called {\em the superextension} of $X$.

The following diagram describes the inclusion relations between subspaces $N_kX$, $N_{<\w}X$, $\Fil(X)$, $\lambda X$ and $\beta X$ of $G(X)$ (an arrow $A\to B$ means that $A$ is a subset of $B$).

\begin{picture}(300,70)(-70,10)
\put(5,60){$\Fil(X)\to N_{<\w}X\to N_{k}X\to N_{2}X\to G(X)$}
\put(13,20){$\beta X$}
\put(20,30){\vector(0,1){25}}
\put(33,22){\vector(1,0){97}}
\put(135,20){$\lambda X$}
\put(140,30){\vector(0,1){25}}
\end{picture}

For a normal space $X$ all the subspaces from this diagram are closed in $G(X)$, see \cite{Gav}.

For a non-compact space $X$ we can also consider the intersections $$
\begin{aligned}
\Fil^\circ(X)=&\Fil(X)\cap G^\circ(X),\quad N^\circ_{<\w}(X)=N_{<\w}(X)\cap G^\circ(X),\\
N_k^\circ (X)=&N_k(X)\cap G^\circ(X),\quad \lambda^\circ_k(X)=\lambda_k(X)\cap G^\circ(X), \mbox{ and}\\
\beta^\circ (X)=&\beta X\cap G^\circ(X)=\beta X\setminus X.
\end{aligned}
$$ Elements of those sets will be called free filters, free centered inclusion hyperspaces, free $k$-linked inclusion hyperspaces, etc. For a normal locally compact space $X$ the subsets 
$\Fil^\circ(X)$, $N^\circ_{<\w}(X)$, $N_k^\circ (X)$, $\lambda^\circ(X)=\lambda^\circ_2(X)$, and $\beta^\circ (X)$ are closed in $G(X)$, see \cite{Gav}. In contrast, $\lambda_k^\circ(\IN)$ is not closed in $G(\IN)$ for $k\ge 3$, see \cite{Gav2}.

\subsection{The inner algebraic structure of $G(X)$} In this subsection we discuss the algebraic
 structure of the space of inclusion hyperspaces $G(X)$ over a topological space $X$. The space of inclusion hyperspaces $G(X)$
 possesses two binary operations $\cup$, $\cap$, and one unary  operation $$\perp:G(X)\to G(X),\; \perp:\F\mapsto\F^\perp=\{E\clin X:\forall F\in\F\;\; E\cap F\ne\emptyset\}$$called the transversality map.
  These three operations are continuous and turn $G(X)$ into a symmetric lattice, see \cite{Gav}.

\begin{defn}
A {\em symmetric lattice} is a complete distributive lattice $(L,\vee ,\wedge)$
endowed with an additional unary operation $\perp :L \to L$, $\perp
:x \mapsto x^\perp$, that is an involutive anti-isomorphism in the
sense that
\begin{itemize}
\item $x^{\perp\perp}=x$ \text{for all} $x \in L$; \item $(x\vee
y)^\perp =x^\perp \wedge y^\perp$; \item $(x\wedge y)^\perp
=x^\perp \vee y^\perp$;
\end{itemize}
\end{defn}

The smallest element of
the lattice $G(X)$ is the inclusion hyperspace $\{X\}$ while the
largest is $\exp(X)$.

 For a discrete space $X$ the set $G(X)$ of all
inclusion hyperspaces on $X$ is a subset of the double power-set
$\mathcal{P}(\mathcal{P}(X))$ (which is a complete distributive
lattice) and is closed under the operations of union and
intersection (of arbitrary families of inclusion hyperspaces).

Since each inclusion hyperspace is a union of filters and each
filter is an intersection of ultrafilters, we obtain the following
proposition showing that the lattice $G(X)$ is a rather natural
object.

\begin{stat}
For a discrete space $X$ the lattice $G(X)$ coincides with the
smallest complete sublattice of $\mathcal{P}(\mathcal{P}(X))$
containing all ultrafilters.
\end{stat}

\section{Extending algebraic operations to inclusion hyperspaces}\label{s2}

In this section, given a binary (associative) operation $*:X\times X\to X$ on a discrete space $X$ we extend this operation to a right-topological (associative) operation on $G(X)$. This can be done in two steps by analogy with the extension of the operation to the Stone-\v Cech compactification $\beta X$ of $X$.

First, for each element $a \in X$ consider
the left shift $L_a:X\to X$, $L_{a}(x)=a\ast x$ and extend it to a continuous
map $\bar{L}_{a}:\beta X\rightarrow\beta X$ between the Stone-\v Cech compactifications of $X$. Next, apply to this extension the functor $G$ to obtain the
continuous map $G\bar{L}_{a}:G(\beta X)\rightarrow G(\beta
X)$. Clearly, for every inclusion hyperspace $\mathcal{F} \in
G(\beta X)$ the inclusion hyperspace $G\bar{L}_{a}(\mathcal{F})$
has a base $\{a\ast F \;|\;F \in \mathcal{F}\}$. Thus, we have
defined the product $a\ast \mathcal{F}= G\bar{L}_{a}(\mathcal{F})$
of the element $a \in X$ and the inclusion hyperspace
$\mathcal{F}$.

Further, for each inclusion hyperspace $\mathcal{F} \in G(\beta
X)=G(X)$ we can consider the map $R_{\mathcal{F}}:X\rightarrow
G(\beta X)$ defined by the formula $R_{\mathcal{F}}(x)=x\ast
\mathcal{F}$ for every $x \in X$. Extend the map $R_{\mathcal{F}}$
to a continuous map $\bar{R}_{\mathcal{F}}: \beta X \rightarrow
G(\beta X)$ and apply to this extension the functor $G$ to obtain a map
$G\bar{R}_{\mathcal{F}}: G(\beta X) \rightarrow G^{2}(\beta X)$.
Finally, compose the map $G\bar R_{\mathcal F}$ with the multiplication $\mu X=\mu
_{G}X:G^2X \rightarrow G(X)$ of the monad $\mathbb G=(G,\eta,\mu)$ and obtain a map $\mu_X\circ G\bar R_{\mathcal F}:G(\beta X)\to G(\beta X)$. For an inclusion hyperspace $\U\in G(\beta X)$, the image
 $\mu _{G}X\circ
G\bar{R}_{\mathcal{F}}(\mathcal{U})$ is called the product of the
inclusion hyperspaces $\mathcal{U}$ and $\mathcal{F}$ and is
denoted by $\mathcal{U}\circ\mathcal{F}$.


It follows from continuity of the maps $G\bar{R}_{\mathcal{F}}$
that the extended binary operation on $G(X)$ is continuous with respect to
the first argument with the second argument fixed. We are going to show that the operation $\circ $ on $G(X)$  nicely agrees with the lattice structure of $G(X)$ and is associative if so is the operation $*$. Also we shall establish an easy formula
$$\U\circ\F=\langle\bigcup_{x\in U}x*F_x:U\in\U,\;\{F_x\}_{x\in U}\subset\F\rangle$$
for calculating the product $\U\circ\F$ of two inclusion hyperspaces $\U,\F$. We start with necessary definitions.

\begin{defn} Let $\star:G(X)\times G(X)\to G(X)$ be a binary operation
on $G(X)$. We shall say that $\star$ {\em respects} the lattice
structure of $G(X)$ if for any $\U,\V,\W\in G(X)$ and $a \in X$
\begin{enumerate}
\item $(\U\cup\V)\star\W=(\U\star\W)\cup(\V\star\W)$; \item
$(\U\cap\V)\star\W=(\U\star\W)\cap(\V\star\W)$; \item
$a\star(\V\cup\W)=(a\star\V)\cup(a\star\W)$; \item
$a\star(\V\cap\W)=(a\star\V)\cap(a\star\W)$.
\end{enumerate}
\end{defn}

\begin{defn}
 We will say that a binary operation $\star:G(X)\times G(X)\to G(X)$ is right-topological if
\begin{itemize}
\item for any $\U\in G(X)$ the right shift $R_\U:G(X)\to G(X)$,
$R_\U:\F\mapsto \F\star\U$, is continuous; \item for any $a\in X$
the left shift $L_a:G(X)\to G(X)$, $L_a:\F\mapsto a\star\F$, is
continuous.
\end{itemize}
\end{defn}

The following uniqueness theorem will be used to find an equivalent description of the induced operation on $G(X)$.

\begin{theo}\label{coincide} Let $\star,\circ:G(X)\times G(X)\to G(X)$ be two right-topological binary
operations that
 respect the lattice structure of $G(X)$. These operations coincide if and only
 if they coincide on the product $X\times X\subset G(X)\times G(X)$.
\end{theo}

\begin{proof}
It is clear that if these operations coincide on $G(X)\times G(X)$,
then they coincide on the product $X\times X$ identified with a subset of $G(X)\times G(X)$. We recall that each point $x\in X$ is identified with the ultrafilter $\gen[ x]$ generated by $x$.

Now assume conversely that $x\star y=x\circ y$ for any two points $x,y\in X\subset G(X)$.
First we check that $a\star \F=a\circ \F$ for any $a \in X$ and $\F \in G(X)$.
Since the left shifts $\F\mapsto a\star\F$ and $\F\mapsto a\circ\F$ are continuous, it suffices to establish the equality $a\star \F=a\circ \F$ for inclusion hyperspaces $\F$ having finite support in $X$ (because the set $G^\bullet(X)$ of all such inclusion hyperspaces is dense in $G(X)$, see \cite{Gav}). Any such a hyperspace $\F$ is generated by a finite family of finite subsets of $X$.

If $\F=\gen[F]$ is generated by a single finite subset $F=\{a_1,\dots,a_n\}\subset X$, then $\F=\bigcap_{i=1}^n \gen[a_i]$ is the finite intersection of principal ultrafilters, and hence $$\gen[a]\star \F=\gen[a]\star\bigcap_{i=1}^n \gen[a_i]=\bigcap_{i=1}^n \gen[a]\star \gen[a_i]=\bigcap_{i=1}^n\gen[a]\circ \gen[a_i]=\gen[a]\circ\bigcap_{i=1}^n\gen[a_i]=\gen[a]\circ \F.$$

If $\F=\gen[F_1,\dots,F_n]$ is generated by finite family of finite sets, then $\F=\bigcup_{i=1}^n\gen[F_i]$ and we can use the preceding case to prove that
$$\gen[a]\star
\F=\gen[a]\star\bigcup_{i=1}^n\gen[F_i]=\bigcup_{i = 1}^n \gen[a]\star{\gen[F_{i}]}=\bigcup_{i = 1}^n \gen[a]\circ\gen[F_{i}]=
\gen[a]\circ\bigcup_{i = 1}^n {\gen[F_{i}]}=
\gen[a] \circ \F.$$

Now fixing any inclusion hyperspace $\U\in G(X)$ by a similar argument one can prove the equality $\F\star\U=\F\circ\U$ for all inclusion hyperspaces $\F\in G^\bullet(X)$ having finite support in $X$. Finally, using  the density of $G^\bullet(X)$ in $G(X)$ and the continuity of right shifts $\F\mapsto \F\circ\U$ and $\F\mapsto \F\star\U$ one can establish the equality $\F\star\U=\F\circ\U$ for all inclusion hyperspaces $\F\in G(X)$.
\end{proof}

The above theorem will be applied to show that the operation $\circ:G(X)\times G(X)\to G(X)$ induced by the  operation $*:X\times X\to X$ coincides with the operation $\star:G(X)\times G(X)\to G(X)$ defined by the formula $$\U \star\V =
\gen[\bigcup\limits_{\substack{
 x \in U}} {x\ast V_x } \colon U\in \U ,\;\{V_x\}_{x\in U} \subset \V ]$$ for  $\U,\V\in G(X)$.

First we establish some properties of the operation $\star$.

\begin{stat}\label{trans1} The operation $\star$ commutes with the transversality operation in the sense that $(\mathcal{U} \star\mathcal{V})^{\bot}=\mathcal{U}^{\bot}
 \star\mathcal{V}^{\bot}$ for any $\mathcal{U} , \mathcal{V}\in G(X)$.
\end{stat}

\begin{proof} To prove that $\U^\bot\star\V^\bot\subset(\U\star\V)^\bot$, take any element $A \in \mathcal{U}^{\bot}
 \star\mathcal{V}^{\bot}$. We should check that $A$ intersects each set $B\in \U\star\V$. Without loss of generality, the sets $A$ and $B$ are of the basic form:
$$A=\bigcup_{\substack{x \in F}} {x\ast G_x}\mbox{ \ for some sets $F \in
\mathcal{U}^{\bot}$ and $\{G_x\}_{x\in F}\subset \mathcal{V}^{\bot}$}$$ and
$$B=\bigcup_{\substack{x \in U}} {x\ast V_x}\mbox{ \ for some sets $U \in
\mathcal{U}$ and $\{V_x\}_{x\in U} \subset \mathcal{V}$}.$$

Since $U\in\U$ and $F\in\U^\bot$, the intersection $F\cap U$ contains some point $x$. For this point $x$ the sets $V_x\in \V$ and $G_x\in\V^\bot$ are well-defined and their intersection $V_x\cap G_x$ contains some point $y$. Then the intersection $A\cap B$ contains the point $x*y$ and hence is not empty, which proves that $A\in (\U\star\V)^\bot$.
\smallskip

 To prove that $(\U\star\V)^\bot\subset \U^\bot\star\V^\bot$, fix a set $A \in (\mathcal{U} \star
\mathcal{V})^{\bot}$. We claim that the set $$F=\{x \in
X: x^{-1} A \in \mathcal{V}^{\bot}\}$$ belongs to $\mathcal{U}^{\bot}$ (here $x^{-1}A=\{y\in X:x\ast y\in A\}$). Assuming conversely that $F \notin
\mathcal{U}^{\bot}$, we would find a set $U \in \mathcal{U}$ with $F\cap U =\emptyset$. By the definition of $F$, for each $x\in U$ the set $x^{-1} A\notin \V^\bot$ and thus we can find a set $V_x\in\V$ with empty intersection $V_x\cap x^{-1}A$. By the definition of the product $\U\star\V$, the set $B=\bigcup_{x\in U}x*V_x$ belongs to $\U\star \V$ and hence intersects the set $A$. Consequently, $x*y\in A$ for some $x\in U$ and $y\in V_x$. The inclusion $x*y\in A$ implies that $y\in x^{-1}A\subset X\setminus V_x$, which is a contradiction proving that $F\in\U^\bot$.
Then the sets $A\supset \bigcup_{x\in F}x\ast x^{-1} A$ belong to $\U^\bot\star\V^\bot$.
 \end{proof}

\begin{stat}\label{RT}
  The equality
$(\mathcal{U} \cap \mathcal{V})\star \mathcal{W}= (\mathcal{U}\star
\mathcal{W})\cap(\mathcal{V}\star \mathcal{W})$ holds for any $\mathcal{U} , \mathcal{V},\mathcal{W}\in G(X)$.
\end{stat}

\begin{proof}
It is easy to show that $(\mathcal{U} \cap \mathcal{V})\star
\mathcal{W}\subset (\mathcal{U}\star \mathcal{W})\cap
(\mathcal{V}\star\mathcal{W})$.

 To prove the reverse inclusion, fix a set $F \in (\mathcal{U}\star
\mathcal{W})\cap (\mathcal{V}\star \mathcal{W})$. Then $$F\supset
\bigcup\limits_{\substack{x \in U}} {x\ast W'_{x}}\mbox{ and }F\supset
\bigcup\limits_{\substack{y \in V}} {y\ast W_{y}''}$$ for some $U \in
\mathcal{U}$, $\{W'_{x}\}_{x\in U}\subset \mathcal{W}$, and  $V \in
\mathcal{V}$, $\{W_{y}''\}_{y\in V}\subset  \mathcal{W}$. Since
$\mathcal{U} , \mathcal{V}$ are inclusion
hyperspaces, $U\cup V \in \mathcal{U} \cap \mathcal{V}$. For each $z\in U\cup V$ let $W_z=W'_z$ if $z\in U$ and $W_z=W''_z$ if $z\notin U$. It follows that
 $F\supset\bigcup\limits_{z\in U\cup V}z\ast W_z$ and hence $F\in (\mathcal{U} \cap \mathcal{V})\star \mathcal{W}$.
\end{proof}

By analogy one can prove

\begin{stat}\label{lattice2} For any $\U,\V,\W\in G(X)$ and $a \in X$
$$a\star(\V\cup\W)=(a\star\V)\cup(a\star\W)\mbox{ \ and \ }
a\star(\V\cap\W)=(a\star\V)\cap(a\star\W).$$
\end{stat}

Combining Propositions~\ref{trans1} and \ref{RT} we get

\begin{corollary}\label{lattice3} For any $\U,\V,\W\in G(X)$ we get
$$(\mathcal{U} \cup \mathcal{V})\star \mathcal{W}=(\mathcal{U}\star \mathcal{W})\cup
(\mathcal{V}\star \mathcal{W}).$$
\end{corollary}

\begin{proof} Indeed,
$$
\begin{aligned}
(\mathcal{U} \cup \mathcal{V})\star \mathcal{W}=\;&\big(((\mathcal{U}
\cup \mathcal{V})\star \mathcal{W})^{\bot}\big)^{\bot}=((\mathcal{U}
\cup \mathcal{V})^{\bot}\star
\mathcal{W}^{\bot})^{\bot}=\\
=&((\mathcal{U}^{\bot} \cap
\mathcal{V}^{\bot})\star
\mathcal{W}^{\bot})^{\bot}=((\mathcal{U}^{\bot}\star
\mathcal{W}^{\bot})\cap (\mathcal{V}^{\bot}\star
\mathcal{W}^{\bot}))^{\bot}=\\
=&(\mathcal{U}^{\bot}\star
\mathcal{W}^{\bot})^{\bot}\cup (\mathcal{V}^{\bot}\star
\mathcal{W}^{\bot})^{\bot}=(\mathcal{U}\star \mathcal{W})\cup
(\mathcal{V}\star \mathcal{W}).
\end{aligned}$$
\end{proof}

\begin{stat}\label{p2.8}
The operation $$\star:G(X)\times G(X)\to G(X),\quad \U \star\V =
\gen[\bigcup_{\substack{
 x \in U}} {x\ast V_x } :U\in \U ,\;\{V_x\}_{x\in U}\subset \V],$$ respects the lattice structure of $G(X)$ and is
right-topological.
\end{stat}

\begin{proof} Propositions~\ref{RT}, \ref{lattice2} and Corollary~\ref{lattice3} imply that the operation $\star$ respects the lattice structure of $G(X)$.

So it remains to check that the operation $\star$ is right-topological.
First we check that for any $\U\in G(X)$ the right shift $R_\U:G(X)\to
G(X)$, $R_\U:\F\mapsto \F\star\U$, is continuous. 

Fix any inclusion hyperspaces $\F,\U\in G(X)$ and let $W^+$ be a sub-basic neighborhood of their product $\F\star\U$. Find sets $F\in\F$ and $\{U_x\}_{x\in F}\subset\U$ such that
$\bigcup\limits_{\substack{
 x \in F}} {x\ast U_x } \subset W$. Then $F^+$ is a neighborhood of $\F$ with
 $F^{+}\star\U \subset W^{+}$.

Now assume that $\F\star\U \in W^{-}$ for some $W\subset X$. Observe that for any inclusion hyperspace $\V\in G(X)$ we get the equivalences $\V\in W^{-} \Leftrightarrow W\in
\V ^{\perp} \Leftrightarrow \V ^{\perp}\in W^{+}$. Consequently, $\F\star\U\in W^-$ is equivalent to
$\F^\perp\star\U^\perp =(\F\star\U)^\perp \in W^{+}$. The preceding case yields
a neighborhood $O(\F^\perp)$ such that
$O(\F^\perp)\star\U^\perp\in W^{+}$. Now the continuity of the transversality operation implies that $O(\F^\perp)^\perp$ is a neighborhood of $\F$ with
$O(\F^\perp)^\perp\star\U\in W^{-}$.

Finally, we prove that  for every $a\in X$ the left shift $L_a:G(X)\to G(X)$,
$L_a:\F\mapsto a\star\F$, is continuous. Given a sub-basic open set $W^+\subset G(X)$ note that $L_a^{-1}(W^+)$ is open because $L_a^{-1}(W^+)=(a^{-1}W)^+$ where $a^{-1}W=\{x\in X:a*x\in W\}$. On the other hand, $a\star \F\in W^-$ is equivalent to $a\star\F^\perp=(a\star \F)^\perp\in (W^-)^\perp=W^+$ which implies that the preimage
$$L_a^{-1}(W^-)=(L_a(W^+))^\perp$$ is also open.
\end{proof}

The operation $\circ$ has the same properties.

\begin{stat}
The operation $\circ:G(X)\times G(X)\to G(X)$, $\U \circ\V =\mu
_{G}X\circ G\bar{R}_{\mathcal{F}}(\mathcal{U})$ respects the
lattice structure of $G(X)$ and is right-topological.
\end{stat}

\begin{proof}
For any $\U \in G(X)$ the right shift $R_{\U}=\mu_{G(X)}\circ G\bar R_{\U}: G(X) \to G(X)$, $R_{\U}:
\F \mapsto \F\circ\U$ is continuous being the composition of continuous maps.
Next for any $a \in X$  and $\F\in G(X)$ we have
$L_{a}(\F)=a\circ\F=\mu_{G}X(\gen[a]\ast\F)=\mu_{G}X(\gen[a\ast\F])=a\ast\F=G\bar{L}_a(\F)$ and the map $L_{a}\equiv G\bar L_a$ is
continuous.

It is known (and easy to verify) that the multiplication $\mu_{G(X)}:G^2(X)\to G(X)$ is a lattice homomorphism in the sense that $$\mu _{G(X)}(\U \cup \V)=\mu _{G(X)}(\U)\cup\mu
_{G}X(\V)\mbox{ \ and \ }\mu _{G(X)}(\U \cap \V)=\mu _{G(X)}(\U)\cap\mu _{G(X)}(\V)$$
for any $\U ,\V \in G(X)$.
Then for any $\U,\V,\W\in G(X)$ and $a \in X$ we get
$$\begin{aligned}
(\U\cup\V)\circ\W=\;&\mu _{G(X)}\circ
G\bar{R}_{\mathcal{W}}(\mathcal{U}\cup \V)=\mu _{G(X)} (
G\bar{R}_{\mathcal{W}}(\mathcal{U})\cup
G\bar{R}_{\mathcal{W}}(\mathcal{V}))=\\
=&\mu _{G(X)}\circ
G\bar{R}_{\mathcal{W}}(\mathcal{U})\cup\mu _{G(X)}\circ
G\bar{R}_{\mathcal{W}}(\mathcal{V})=(\U\circ\W)\cup(\U\circ\W)
\end{aligned}$$
and similarly $
(\U\cap\V)\circ\W=(\U\circ\W)\cap(\U\circ\W).$

Note that for any $a \in X$ $$a\circ \W =\mu _{G(X)} (
G\bar{R}_{\mathcal{W}}(\gen[a])=\gen[\bar{R}_{\mathcal{W}}(
\{a\})]=\gen[\bar{R}_{\mathcal{W}}( a)]=a\ast\W.$$

Consequently,
$$a\circ(\V\cup\W)=a\ast(\V\cup\W)=
(a\ast\V)\cup(a\ast\W)=(a\circ\V)\cup(a\circ\W)
$$ and similarly $a\circ(\V\cap\W)=(a\circ\V)\cap(a\circ\W)$.
\end{proof}

Since both operations $\circ$ and $\star$ are right-topological and respect the lattice structure of $G(X)$ we may apply Theorem~\ref{coincide} to get

\begin{cons}\label{apparent} For any binary operation $*:X\times X\to X$ the operations
 $\circ$ and $\star$ on $G(X)$ coincide. Consequently, for any inclusion hyperspaces $\U,\V\in G(X)$ their product $\U\circ\V$ is the inclusion hyperspace $$\langle\bigcup_{x\in U}x*V_x:U\in\U,\;\{V_x\}_{x\in U}\subset\V\rangle=\big\{A\subset X:\{x\in X:x^{-1}A\in\V\}\in\U\big\}.$$
\end{cons}

Having the apparent description of the operation $\circ$ we can establish its associativity.

\begin{stat}
If the operation $\ast$ on $X$ is associative, then so is the
induced operation $\circ$ on $G(X)$.
\end{stat}

\begin{proof}
It is necessary to show that $(\mathcal{U}
\circ\mathcal{V})\circ\mathcal{W} =\mathcal{U} \circ
(\mathcal{V}\circ\mathcal{W})$ for any inclusion hyperspaces
$\mathcal{U} ,\mathcal{V} ,\mathcal{W}$.

Take any subset $A \in (\mathcal{U}
\circ\mathcal{V})\circ\mathcal{W}$ and choose a set $B \in
\mathcal{U} \circ\mathcal{V}$ such that $A \supset
\bigcup\limits_{\substack{
 z \in B}} {z\ast W_z }$ for some family $\{W_z\}_{z\in B} \subset \mathcal{W}$.
 Next, for the set $B \in \mathcal{U}
\circ\mathcal{V}$ choose a set $U \in \mathcal{U}$ such that $B
\supset \bigcup\limits_{\substack{
 x \in U}} {x\ast V_x }$ for some family $\{V_x\}_{x\in U} \subset\mathcal{V}$.
 It is clear that for each $x \in U$ and $y \in V_x$ the product
 $x\ast y$ is in $B$ and hence $W_{x\ast y}$ is defined. Consequently,
 $\bigcup\limits_{\substack{
 y \in V_x}} {y\ast W_{x\ast y} } \in \mathcal{V}\circ\mathcal{W}$
 for all $x \in U$ and hence $\bigcup\limits_{\substack{
 x \in U}} {x}\ast (\bigcup\limits_{\substack{
 y \in V_x}} {y\ast W_{x\ast y} }) \in \mathcal{U} \circ
(\mathcal{V}\circ\mathcal{W})$. Since $\bigcup\limits_{\substack{
 x \in U}} \bigcup\limits_{\substack{
 y \in V_x}} {x\ast y\ast W_{x\ast y} } \subset A$, we get $A \in \mathcal{U} \circ
(\mathcal{V}\circ\mathcal{W})$. This proves the inclusion
$(\mathcal{U} \circ\mathcal{V})\circ\mathcal{W} \subset
\mathcal{U} \circ (\mathcal{V}\circ\mathcal{W})$.

To prove the reverse inclusion, fix a set $A \in \mathcal{U} \circ
(\mathcal{V}\circ\mathcal{W})$ and choose a set $U \in \mathcal{U}$
such that $A \supset \bigcup\limits_{\substack{
 x \in U}} {x\ast B_x }$ for some family $\{B_x\}_{x\in U} \subset \mathcal{V}\circ\mathcal{W}$. Next, for each $x \in U$ find a set $V_x \in \mathcal{V}$
 such that $B_x \supset \bigcup\limits_{\substack{
 y \in V_{x}}} {y\ast W_{x, y}}$ for some family $\{W_{x, y}\}_{y\in V_x} \subset \mathcal{W}$. Let $Z=\bigcup\limits_{\substack{
 x \in U}} {x\ast V_x }$. For each $z \in Z$ we can find $x \in U$ and
 $y \in V_x$ such that $z=x\ast y$ and put $W_z = W_{x, y}$. Then $Z \in \mathcal{U}
\circ\mathcal{V}$ and $\bigcup\limits_{\substack{
 z \in Z}} {z\ast W_z } \in (\mathcal{U}
\circ\mathcal{V})\circ\mathcal{W}$. Taking into account
$\bigcup\limits_{\substack{
 z \in Z}} {z\ast W_z } \subset \bigcup\limits_{\substack{
 x \in U}}\bigcup\limits_{\substack{
 y \in V_x}} {x\ast y\ast W_{x ,y }} \subset A$, we conclude $A \in (\mathcal{U}
\circ\mathcal{V})\circ\mathcal{W}$.
\end{proof}

\section{Homomorphisms of semigroups of inclusion hyperspaces}

Let us observe that our construction of extension of a binary operation for $X$ to $G(X)$ works well both for associative and non-associative operations. Let us recall that a set $S$ endowed with a binary operation $*:X\times X \to X$ is called a {\em groupoid}. If the operation is associative, then $X$ is called a {\em semigroup}. In the preceding section we have shown that for each groupoid (semigroup) $X$ the space $G(X)$ is a groupoid (semigroup) with respect to the extended operation.

A map $h:X_1\to X_2$ between two groupoids $(X_1,*_1)$ and $(X_2,*_2)$ is called a {\em homomorphism} if $h(x*_1y)=h(x)*_2h(y)$ for all $x,y\in X_1$.

\begin{stat} For any homomorphism $h:X_1\to X_2$ between groupoids $(X_1,*_1)$ and $(X_2,*_2)$ the induced map $Gh:
G(X_1)\rightarrow G(X_2)$ is a homomorphism of the groupoids $G(X_1)$, $G(X_2)$.
\end{stat}

\begin{proof}
Given two inclusion hyperspaces $\mathcal{U} ,\mathcal{V}\in G(X_1)$ observe that
$$
\begin{aligned}
Gh(\mathcal{U}\circ_{1}\mathcal{V})=\;&
Gh(\langle\bigcup\limits_{\substack{x \in U}} {x\ast_1 V_{x}}: U \in
\mathcal{U},\;\{V_{x}\}_{x\in U}\subset\mathcal{V}\rangle)=\\
=\;&\langle h(\bigcup\limits_{\substack{x
\in U}} {x\ast_1 V_{x}}):U \in \mathcal{U},\;\{V_{x}\}_{x\in U}\subset\mathcal{V}\rangle)=\\
=\;&\langle\bigcup\limits_{\substack{x \in U}} {h(x)\ast_2
h(V_{x})}:U \in \mathcal{U},\;\{V_{x}\}_{x\in U}\subset \mathcal{V}\rangle=\\
=\;&\langle\bigcup\limits_{\substack{x \in h(U)}} {x\ast_2
h(V_{x})}:U
\in \mathcal{U},\; \{h(V_{x})\}_{x\in U}\subset Gh(\mathcal{V})\rangle=\\
=\;&\langle h(U):U \in
\mathcal{U}\rangle \circ_2 \langle h(V):V \in \mathcal{V}\rangle=Gh(\mathcal{U})\circ_{2}Gh(\mathcal{V}).
\end{aligned}
$$
\end{proof}

Reformulating Proposition~\ref{trans1} in terms of homomorphisms, we obtain

\begin{stat}\label{subtrans} For any groupoid $X$ the transversality map $\perp:G(X)\to G(X)$ is a homomorphism of the groupoid $G(X)$.
\end{stat}

\section{Subgroupoids of $G(X)$}

In this section we shall show that for a groupoid $X$ endowed with the discrete topology all (topologically) closed subspaces of $G(X)$ introduced in Section~\ref{subspaces} are subgroupoids of $G(X)$. A subset $A$ of a groupoid $(X,*)$ is called a {\em subgroupoid} of $X$ if $A*A\subset A$, where $A*A=\{a*b:a,b\in A\}$.

We assume that $*:X\times X\to X$ is a binary operation on a discrete space $X$ and $\circ:G(X)\times G(X)\to G(X)$ is the extension of $*$ to $G(X)$. Applying Proposition~\ref{subtrans} we obtain

\begin{stat} If $S$ is a subgroupoid of $G(X)$, then $S^\perp$ is a subgroupoid of $G(X)$ too.
\end{stat}

Our next propositions can be easily derived from Corollary~\ref{apparent}.

\begin{stat} The sets $\Fil(X)$, $N_{<\w}(X)$ and $N_k(X)$, $k\ge 2$, are  subgroupoids in $G(X)$.
\end{stat}




\begin{stat} The Stone-\v Cech extension $\beta X$ and the superextension $\lambda X$ both are closed subgroupoids in $G(X)$.
\end{stat}

\begin{proof} The superextension $\lambda X$ is a subgroupoid of $G(X)$ being the intersection $\lambda(X)=N_2(X)\cap(N_2(X))^\perp$ of two subgroupoids of $G(X)$. By analogy, $\beta X=\Fil(X)\cap \lambda(X)$ is a subgroupoid of $G(X)$.
\end{proof}

\begin{rem} In contrast to $\lambda X$ for $k\ge 3$ the subset $\lambda_k(X)$ need not be a subgroupoid of $G(X)$. For example, for the cyclic group $\IZ_5=\{0,1,2,3,4\}$ the subset $\lambda_3(\IZ_5)$  of $G(\IZ_5)$ contains a maximal 3-linked system 
$$\mathcal L=\langle \{0,1,2\},\{0,1,4\}, \{0,2,4\},\{1,2,4\}\rangle$$ whose square $$\mathcal L*\mathcal L=\langle \{1,2,4,5\},\{0,2,3,4\}, \{0,1,3,4\},\{0,1,2,4\},\{0,1,2,3\}\rangle$$
is not maximal 3-linked.
\end{rem} 

By a direct application of Corollary~\ref{apparent} we can also prove

\begin{stat} The set $G^\bullet(X)$ of all inclusion hyperspaces with finite support is a subgroupoid in $G(X)$.
\end{stat} 

Finally we find conditions on the operation $*$ guaranteeing that the subset $G^\circ(X)$ of free inclusion hyperspaces is a subgroupoid of $G(X)$.

\begin{stat}\label{Gcirc} Assume that for each $b\in X$ there is a finite subset $F\subset X$ such that for each $a\in X\setminus F$ the set $a^{-1}b=\{x\in X:a*x=b\}$ is finite. Then the set $G^\circ(X)$ is a closed subgroupoid in $G(X)$ and consequently, $\Fil^\circ(X)$, $\lambda^\circ(X)$, $\beta^\circ(X)$ all are closed subgroupoids in $G(X)$.
\end{stat}

\begin{proof} Take two free inclusion hyperspaces $\A,\mathcal B\in G(X)$ and a subset $C\in \A\circ\mathcal B$. We should prove that $C\setminus K\in\A\circ\mathcal B$ for each compact subset $K\subset X$. Without loss of generality, the set $C$ is of basic form: $C=\bigcup_{a\in A}a*B_a$ for some set $A\in \A$ and some family $\{B_a\}_{a\in A}\subset\mathcal B$.

Since $X$ is discrete, the set $K$ is finite. It follows from our assumption that there is a finite set $F\subset X$ such that for every $a\in X\setminus F$ the set $a^{-1}K=\{x\in X:a*x\in K\}$ is finite. The hyperspace $\A$, being free, contains the set $A'=A\setminus F$.  By the same reason, for each $a\in A'$ the hyperspace $\mathcal B$ contains the set $B'_a=B_a\setminus a^{-1}K$.
Since $C\setminus K\supset \bigcup_{a\in A'}a*B'_a \in\mathcal A\circ\mathcal B$, we conclude that $C\setminus K\in\mathcal A\circ\mathcal B$.
\end{proof}

\begin{remark} If $X$ is a semigroup, then $G(X)$ is a semigroup and all the subgroupoids considered above are closed subsemigroups in $G(X)$. Some of them are well-known in Semigroup Theory.
In particular, so is the semigroup $\beta X$ of ultrafilter and $\beta^\circ(X)=\beta X\setminus X$ of free ultrafilters. The semigroup $\Fil(X)$ contains an isomorphic copy of the global semigroup of $X$, which is the hyperspace $\exp(X)$ endowed with the semigroup operation $A*B=\{a*b:a\in A,\; b\in B\}$.
\end{remark}

\section{Ideals and zeros in $G(X)$}

 A non-empty subset $I$ of a groupoid $(X,*)$ is called an {\em ideal} (resp. {\em right ideal},
 {\em left ideal}\/) if $I*X\cup X*I\subset I$ (resp. $I*X\subset I$, $X*I\subset I$).
An element $O$ of a groupoid $(X,*)$ is called a {\em zero} (resp. {\em left zero}, {\em right zero})
in $X$ if $\{O\}$ is an ideal (resp. right ideal, left ideal) in $X$.
Each right or left zero $z\in X$ is an {\em idempotent} in the sense that $z*z=z$.

For a groupoid $(X,*)$ right zeros in $G(X)$ admit a simple description. We define an inclusion
hyperspace $\mathcal A\in G(X)$ to be {\em shift-invariant} if for every $A\in\mathcal A$ and
$x\in X$ the sets $x*A$ and $x^{-1}A=\{y\in X:x*y\in A\}$ belong to $\mathcal  A$.

\begin{stat}\label{zeros} An inclusion hyperspace $\mathcal A\in G(X)$ is a right zero
in $G(X)$ if and only if $\A$ is shift-invariant.
\end{stat}

\begin{proof}
Assuming that an inclusion hyperspace $\mathcal A\in G(X)$ is
shift-invariant, we shall show that $\mathcal B\circ\A=\A$ for every $\mathcal B\in G(X)$. Take any set $F \in \mathcal
B\circ\A$ and find a set $B\in\mathcal B$ and a family $\{A_x\}_{x\in B}\subset \A$ such that $ \bigcup\limits_{\substack{x \in
B}} {x\ast A_{x}}\subset F$. Since
$\mathcal A\in G(X)$ is shift-invariant,
$\bigcup\limits_{\substack{x \in B}} {x\ast A_{x}} \in \A$ and
thus $F \in \A$. This proves the inclusion $\mathcal B\circ\A\subset\A$. On the other hand, for every $F\in \A$ and every $x\in X$ we get $x^{-1}F \in \A$ and
thus $F\supset\bigcup\limits_{\substack{x \in X}} {x\ast
x^{-1}F} \in \mathcal B \circ \A$. This shows that $\A$ is a right zero of the semigroup $G(X)$.

Now assume that $\A$ is a right zero of $G(X)$. Observe that for every $x\in X$ the equality $\gen[x]\circ\A=\A$ implies $x*A\in\A$ for every $A\in A$.

One the other hand, the equality $\{X\}\circ \A=\A$ implies that for every $A\in\A$ there is a family $\{A_x\}_{x\in X}\subset\A$ such that $\bigcup_{x\in X}x*A_x\subset A$. Then for every $x\in X$ the set $x^{-1}A=\{z\in X:x*z\in A\}\supset A_x\in\A$ belongs to $\A$ witnessing that $\A$ is shift-invariant. 
\end{proof}

By $\invG(X)$ we denote the set of shift-invariant inclusion hyperspaces in $G(X)$.
 Proposition~\ref{zeros} implies that $\A\circ\mathcal B=\mathcal B$ for every
  $\A,\mathcal B\in \invG(X)$. This means that $\invG(X)$ is a rectangular semigroup.

We recall that a semigroup $(S,*)$ is called  {\em rectangular}
(or else a {\em semigroup of right zeros}) if $x*y=y$ for all $x,y\in S$.

\begin{stat}\label{idealfull} The set $\invG(X)$ is closed in $G(X)$, is a rectangular subsemigroup
of the groupoid $G(X)$ and is closed complete sublattice of the lattice $G(X)$ invariant
under the transversality map. Moreover, if $\invG(X)$ is non-empty, then it is a left
ideal that lies in each right ideal of $G(X)$.
\end{stat}

\begin{proof}
If $\A\in G(X)\setminus\invG(X)$, then there exists $x \in X$ and
$A\in\A$ such that $x*A \notin \A$ or $x^{-1}A\notin\A$. Then $$O(\A)=\{\A'\in G(X): A\in\A' \mbox{  and }(x*A\notin\A'\mbox{  or }x^{-1}A\notin\A)\}$$ is an open neighborhood of $\A$ missing the set $\invG(X)$ and witnessing that the set $\invG(X)$ is closed in $G(X)$.  

Since $\A\circ\mathcal B=\mathcal B$ for every $\A,\mathcal B\in \invG(X)$, the set $\invG(X)$ is a rectangular subsemigroup of the groupoid $G(X)$.

To show that $\invG(X)$ is invariant under the transversality operation, note that for every $\A\in G(X)$ and $\mathcal Z\in\invG(X)$ we get $\A\circ \mathcal Z^\perp=(\A^\perp\circ\mathcal Z)^\perp=\mathcal Z^\perp$ which means that $\mathcal Z^\perp$ is a right zero in $G(X)$ and thus belongs to $\invG(X)$ according to Proposition~\ref{zeros}.

To show that $\invG(X)$ is a complete sublattice of $G(X)$ it is necessary to check that $\invG(X)$ is closed under arbitrary unions and intersections. It is trivial to check that arbitrary union of shift-invariant inclusion hyperspaces is shift-invariant, which means that $\bigcup_{\alpha\in A}\mathcal Z_\alpha\in\invG(X)$ for any family $\{\mathcal Z_\alpha\}_{\alpha\in A}\subset\invG(X)$. Since $\invG(X)$ is closed under the transversality operation we also get $$\bigcap_{\alpha\in A}\mathcal Z_\alpha=\big(\bigcup_{\alpha\in A}\mathcal Z_\alpha^\perp)^\perp\in\invG(X)^\perp=\invG(X).$$

If $\invG(X)$ is not empty, then it is a left ideal in $G(X)$ because it consists of right zeros. Now take any right ideal $I$ in $G(X)$ and fix any element $\mathcal R\in I$. Then for every $\mathcal Z\in\invG(X)$ we get $\mathcal Z=\mathcal R\circ\mathcal Z\in I$ which yields $\invG(X)\subset I$.
\end{proof}

\begin{stat}\label{min} If $X$ is a semigroup and
 $\invG(X)$ is not empty, then $\invG(X)$ is the minimal ideal of $G(X)$.
\end{stat}

\begin{proof} In light of the preceding proposition, it suffices to check that $\invG(X)$ is a right ideal. Take any inclusion hyperspaces $\A\in\invG(X)$ and $\mathcal B\in G(X)$ and take any set $F\in\mathcal A\circ\mathcal B$. We need to show that the sets $x*F$ and $x^{-1}F$ belong to $\mathcal A\circ\mathcal B$. Without loss of generality, $F$ is of the basic form:
$$F=\bigcup_{a\in A}a*B_a$$ for some set $A\in\A$ and some family $\{B_a\}_{a\in A}\subset\mathcal B$. The associativity of the semigroup operation on $S$ implies that
$$x*F=\bigcup_{a\in A}x*a*B_a=\bigcup_{z\in x*A}z*B_{a(z)}\in\A\circ\mathcal B$$where $a(z)\in \{a\in A: x*a=z\}$ for $z\in x*A$. To see that $x^{-1}F\in\mathcal A$ observe that the set $A'=\bigcup_{z\in x^{-1}A}z*B_{xz}$ belongs to $\A$ and  each point $a'\in A'$ belongs to the set $z*B_{xz}$ for some $z\in x^{-1}A$.
Then $x*a'\in x*z*B_{xz}\subset F$ and hence $\A\ni A'\subset x^{-1}F$, which yields the desired inclusion $x^{-1}F\in\A$.\end{proof}

Now we find conditions on the binary operation $*:X\times X\to X$ guaranteeing
that the set $\overset{\leftrightarrow}{G}(X)$ is not empty.
By $\min GX=\{X\}$ and $\max GX=\{A\subset X:A\ne\emptyset\}$ we denote the minimal
 and maximal elements of the lattice $G(X)$.

\begin{stat}\label{maxmin} For a groupoid $(X,*)$ the following conditions are equivalent:
\begin{enumerate}
\item $\min GX\in\invG(X)$;
\item $\max GX\in\invG(X)$;
\item for each $a,b\in X$ the equation $a*x=b$ has a solution $x\in X$.
\end{enumerate}
\end{stat}

\begin{proof}
$(1) \Rightarrow (3)$ Assuming that $\min GX\in\invG(X)$ and applying Proposition~\ref{zeros} observe that for every $a\in X$ the equation $\gen[a]\circ\{X\}=\{X\}$ implies that for every $b\in X$ the equation $a*x=b$ has a solution.
\smallskip

$(3) \Rightarrow (1)$ If for every $a,b\in X$ the equation $a*x=b$ has a solution, then $a*X=X$ and hence $\F\circ\{X\}=\{X\}$ for all $\F\in G(X)$.
This means that $\{X\}=\min G(X)$ is a right zero in $G(X)$ and hence belongs to $\invG(X)$ according to Proposition~\ref{zeros}.
\smallskip

$(2)\Rightarrow (3)$ Assume that $\max G(X)\in\invG(X)$ and take any points $a,b\in X$. 
Since $\gen[a]\circ\max G(X)=\max G(X)\ni\{b\}$, there is a non-empty set $X_a\in\max G(X)$ with $a*X_a\subset\{b\}$. Then any $x\in X_a$ is a solution of $a*x=b$.
\smallskip

$(3)\Rightarrow (2)$ Assume that for every $a,b\in X$ the equation $a*x=b$ has a solution. To show that $\F\circ\max G(X)=\max G(X)$ it suffices to check that $\max G(X)\subset\F\circ \max G(X)$. Take any set $B\in\max G(X)$ and any set $F\in\F$.
For every $a\in F$ find a point $x_a\in X$ with $a*x_a\in B$. Then the sets $\bigcup_{a\in F}a*\{x_a\}\subset B$ belong to $\F\circ \max G(X)$, which yields the desired inclusion  $\max G(X)\subset\F\circ \max G(X)$.
\end{proof}

By analogy we can establish a similar description of zeros and the minimal ideal in the semigroup $G^\circ(X)$ of free inclusion hyperspaces.

\begin{stat}\label{ideal} Assume that $(X,*)$ is an infinite groupoid such that for each $b\in X$ there is a finite subset $F\subset X$ such that for each $a\in X\setminus F$ the set $a^{-1}b=\{x\in X:a*x=b\}$ is finite and not empty. Then
\begin{enumerate}
\item $G^\circ(X)$ is a closed subgroupoid of $G(X)$;
\item $G^\circ(X)$ is a left ideal in $G(X)$ provided if for each $a,b\in X$ the set $a^{-1}b$ is finite;
\item the set $\inv[G^\circ](X)=\invG(X)\cap G^\circ(X)$ of shift-invariant free inclusion hyperspaces is the minimal ideal in $G^\circ(X)$;
\item the set $\inv[G^\circ](X)$ is a rectangular subsemigroup of the groupoid $G(X)$ and is closed complete sublattice of the lattice $G(X)$ invariant under the transversality map.
\end{enumerate}
\end{stat}

\begin{remark} It follows from Propositions~\ref{idealfull} and \ref{ideal} that the minimal ideals of the semigroups $G(\IZ)$ and $G^\circ(X)$ are closed. In contrast, the minimal ideals of the semigroups $\beta\IZ$ and $\beta^\circ\IZ=\beta\IZ\setminus\IZ$ are not closed, see \cite[\S 4.4]{HS}.
\end{remark}

Minimal left ideals of the semigroup $\beta^\circ(\IZ)$ play an important role in Combinatorics of Numbers, see \cite{HS}. We believe that the same will happen for the semigroup $\lambda^\circ(\IZ)$. The following proposition implies that minimal left ideals of $\lambda^\circ(\IZ)$ contain no ultrafilter!

\begin{stat}\label{ideall} If a groupoid $X$ admits a homomorphism $h:X\to\IZ_3$ such that for every $y\in\IZ_3$ the preimage $h^{-1}(y)$ is not empty (is infinite) then each minimal left ideal $I$ of $\lambda(X)$ (of $\lambda^\circ(X)$) is disjoint from $\beta(X)$ .
\end{stat}

\begin{proof} It follows that the induced map $\lambda h:\lambda(X)\to\lambda(\IZ_3)$
is a surjective homomorphism. Consequently, $\lambda h(I)$ is a minimal left ideal in $\lambda(\IZ_3)$. Now observe that $\lambda(\IZ_3)$ consists of four maximal linked inclusion hyperspaces. Besides three ultrafilters there is a maximal linked inclusion hyperspace $\mathcal L_\vartriangle= \langle\{0,1\},\{0,2\},\{1,2\}\rangle$ where $\IZ_3=\{0,1,2\}$. One can check that $\{\mathcal L_\vartriangle\}$ is a zero of the semigroup $\lambda(\IZ_3)$. Consequently, $\lambda(h)(I)=\{\mathcal L_\vartriangle\}$, which implies that $I\cap\beta(X)=\emptyset$.

Now assume that for every $y\in\IZ_3$ the preimage $h^{-1}(y)$ is infinite. We claim that $\lambda h(\lambda^\circ(X))=\lambda(\IZ_3)$. Take any maximal linked inclusion hyperspace $\mathcal L\in\lambda(\IZ_3)$. If $\mathcal L$ is an ultrafilter supported by a point $y\in\IZ_3$, then we can take any free ultrafilter $\U$ on $X$ containing the infinite set $h^{-1}(y)$ and observe that $\lambda h(\U)=\mathcal L$. It remains to consider the case $\mathcal L=\mathcal L_\vartriangle$. Fix free ultrafilters $\U_0,\U_1,\U_2$ on $X$ containing the sets $h^{-1}(0)$, $h^{-1}(1)$, $h^{-1}(2)$, respectively. Then $\mathcal L=(\U_0\cap\U_1)\cup(\U_0\cap\U_2)\cup(\U_1\cap\U_2)$ is a free maximal linked inclusion hyperspace whose image $\lambda h(\mathcal L_X)=\mathcal L_\vartriangle$.

Given any minimal left ideal $I\subset \lambda^\circ(X)$ we obtain that the image $\lambda h(I)$, being a minimal left ideal of $\lambda(\IZ_3)$ coincides with $\{\mathcal L_\triangle\}$ and is disjoint from $\beta(\IZ_3)$. Consequently, $I$ is disjoint from $\beta(X)$.
\end{proof}

\section{The center of $G(X)$}

In this section we describe the structure of the center  of the groupoid $G(X)$ for each (quasi)group $X$. By definition, the {\em center} of a groupoid $X$ is the set $$C=\{x\in X:\forall y\in X\;\; xy=yx\}.$$

A groupoid $X$ is called a {\em quasigroup} if for every $a,b\in X$ the system of equations $a*x=b$ and $y*a=b$ has a unique solution $(x,y)\in X\times X$. It is clear that each group is a quasigroup. On the other hand, there are many examples of quasigroups, not isomorphic to groups, see \cite{Qg1}, \cite{Qg2}.

\begin{theo}  Let $X$ be a quasigroup. If an inclusion hyperspace $\C\in G(X)$ commutes with the extremal elements $\max G(X)$ and $\min G(X)$ of $G(X)$, then $\C$ is a principal ultrafilter.
\end{theo}

\begin{proof} By Proposition~\ref{maxmin}, the inclusion hyperspaces $\max G(X)$ and $\min G(X)$ are right zeros in $G(X)$ and thus $\max G(X)\circ\C=\C\circ \max G(X)=\max G(X)$ and $\min G(X)\circ \C=\C\circ\min G(X)=\min G(X)$. It follows that for every $b\in X$ we get $\{b\}\in \max G(X)=\max G(X)\circ\C$, which means that $a*C\subset\{b\}$ for some $C\in\C$ and some $a\in X$. Since the equation $a*y=b$ has a unique solution $y\in X$, the set $C$ is a singleton, say $C=\{c\}$. It remains to prove that $\C$ coincides with the principal ultrafilter $\langle c\rangle$ generated by $c$. Assuming the converse, we would conclude that $X\setminus\{c\}\in\C$. By our hypothesis, the equation $y*c=c$ has a unique solution $y_0\in X$. Since the equation $y_0*x=c$ has a unique solution $x=c$, $y_0*(X\setminus\{c\})\subset X\setminus\{c\}$. Letting $C_x=\{c\}$ for all $x\in X\setminus \{y_0\}$ and $C_{x}=X\setminus\{c\}$ for $x=y_0$, we conclude that $X\setminus\{c\}\supset\bigcup_{x\in X}x*C_x\in\min G(X)\circ\C=\C\circ\min G(X)=\min G(X)$, which is not possible.
\end{proof}

\begin{corollary} For any quasigroup $X$ the center of the groupoid $G(X)$ coincides with the center of $X$.
\end{corollary}

\begin{proof} If an inclusion hyperspace $\C$ belongs to the center of the groupoid $G(X)$, then $\C$ is a principal ultrafilter generated by some point $c\in X$. Since $\C$ commutes with all the principal ultrafilters, $c$ commutes with all elements of $X$ and thus $c$ belongs to the center of $X$. 

Conversely, if $c\in X$ belongs to the center of $X$, then for every inclusion hyperspace $\F\in G(X)$ we get
$$c\circ\F=\{c*F:F\in\F\}=\{F*c:F\in\F\}=\F\circ c,$$
which means that (the principal ultrafilter generated by) $c$ belongs to the center of the groupoid $G(X)$.
\end{proof}

\begin{rem} It is interesting to note that for any group $X$ the center of the semigroup $\beta X$ also coincides with the center of the group $X$, see Theorem 6.54 of \cite{HS}. 
\end{rem}

\begin{problem} Given a group $X$ describe the centers of the subsemigroups $\lambda(X)$, $\Fil(X)$, $N_{<\w}(X)$, $N_k(X)$, $k\ge 2$ of the semigroup $G(X)$. 
Is it true that the center of any subsemigroup $S\subset G(X)$ with $\beta(X)\subset S=S^\perp$ coincides with the center of $X$?
\end{problem}

\begin{rem} Let us note that the requirement $S=S^\perp$ in the preceding question is essential: for any nontrivial group $X$ the center of the (non-symmetric) subsemigroup $X\cup \max G(X)$ of $G(X)$ contains $\max G(X)$ and hence is strinctly larger than the center of the group $X$.
\end{rem}

\begin{problem} Given an infinite group $X$ describe the centers of the semigroups $G^\circ(X)$, $\lambda^\circ(X)$, $\Fil^\circ(X)$, $N^\circ_{<\w}(X)$, and $N^\circ_k(X)$, $k\ge 2$. {\rm  (By Theorem 6.54 of \cite{HS}, the center of the semigroup of free ultrafilters $\beta^\circ(X)$ is empty).}
\end{problem}

\section{The topological center of $G(X)$}

In this section we describe the topological center of $G(X)$. 
By the {\em topological center} of a groupoid $X$ endowed with a topology we understand the set $\Lambda(X)$ consisting of all points $x\in X$ such that the left and right shifts
$$l_x:X\to X, \;\; l_x:z\mapsto xz, \mbox{ \  and  \ }r_x:X\to X, \;\; r_x:z\mapsto zx$$
both are continuous. 

Since all right shifts on $G(X)$ are continuous, the topological center of the groupoid $G(X)$ consists of all inclusion hyperspaces $\F$ with continuous left shifts $l_\F$. 

We recall that $G^\bullet(X)$ stands for the set of inclusion hyperpsaces with finite support. 

\begin{theo}\label{topcenter} For a quasigroup $X$ the topological center of the groupoid $G(X)$ coincides with $G^\bullet(X)$.
\end{theo}

\begin{proof} By Proposition~\ref{p2.8}, the topological center $\Lambda(GX)$ of $G(X)$ contains all principal ultrafilters and is a sublattice of $G(X)$. Consequently, $\Lambda(GX)$ contains the sublatttice $G^\bullet(X)$ of $G(X)$ generated by $X$.

Next, we show that each inclusion hyperspace $\F\in\Lambda(GX)$ has finite support and hence belongs to $G^\bullet(X)$. By Theorem~9.1 of \cite{Gav}, this will follow as soon as we check that both $\F$ and $\F^\perp$ have bases consisting of finite sets. 

Take any set $F\in\F$, choose any point $e\in X$, and consider the inclusion hyperspace $\U=\{U\subset X:e\in F*U\}$. Since for every $f\in F$ the equation $f*u=e$ has a solution in $X$, we conclude that $\{e\}\in\F\circ\U$ and by the continuity of the left shift $l_\F$, there is an open neighborhood $\mathcal O(\U)$ of $\U$ such that 
$\{e\}\in\F\circ\mathcal A$ for all $\mathcal A\in\mathcal O(\U)$.  Without loss of generality, the neighborhood  $\mathcal O(\U)$ is of basic form $$\mathcal O(\U)=U_1^+\cap \dots\cap U_n^+\cap V_1^-\cap \dots\cap V_m^-$$ for some sets $U_1,\dots,U_n\in\U$ and $V_1,\dots,V_m\in\U^\perp$. Take any finite set $A\subset F^{-1}e=\{x\in X:e\in F*x\}$ intersecting each set $U_i$, $i\le n$, and consider the inclusion hyperspace $\A=\langle A\rangle^\perp$. It is clear that $\A\subset U_1^+\cap\dots\cap U_n^+$. 
Since each set $V_j$, $j\le m$, contains the set $F^{-1}e\supset A$, we get also that $\A\in V_1^-\cap\dots\cap V_m^-$. Then $\F\circ \A\ni\{e\}$ and hence there is a set $E\in \F$ and a family $\{A_x\}_{x\in E}\subset\mathcal A$ with $\bigcup_{x\in E}x*A_x\subset\{e\}$. It follows that the set $E\subset eA^{-1}=\{x\in X:\exists a\in A$ with $xa=e\}$ is finite. We claim that $E\subset F$. Indeed, take any point $x\in E$ and find a point $a\in A$ with $x*a=e$. Since $A\subset F^{-1}e$, there is a point $y\in F$ with $e=y*a$. Hence $xa=ya$ and the right cancellativity of $X$ yields $x=y\in F$.
Therefore, using the continuity of the left shift $l_\F$, for every $F\in\F$ we have found a finite subset $E\in\F$ with $E\subset F$. This means that $\F$ has a base of finite sets.

The continuity of the left shift $l_\F$ and Proposition~\ref{trans1} imply the continuity of the left shift $l_{\F^\perp}$. Repeating the preceding argument, we can prove that the inclusion hyperspace $\F^\perp$ has a base of finite sets too. Finally, applying Theorem~9.1 of \cite{Gav}, we conclude that $\F\in G^\bullet(X)$.
\end{proof}

\begin{problem} Given an infinite group $G$ describe the topological center of the subsemigroups $\lambda(X)$, $\Fil(X)$, $N_{<\w}(X)$, $N_k(X)$, $k\ge 2$, of the semigroup $G(X)$. 
Is it true that the topological center of any subsemigroup $S\subset G(X)$ containing $\beta(X)$ coincides with $S\cap G^\bullet(X)$? {\rm (This is true for the subsemigroups $S=G(X)$ (see Theorem~\ref{topcenter})  and $S=\beta(X)$, see Theorems 4.24 and 6.54 of \cite{HS}).}
\end{problem}

\begin{problem} Given an infinite group $X$ describe the topological centers of the semigroups $G^\circ(X)$, $\lambda^\circ(X)$, $\Fil^\circ(X)$, $N^\circ_{<\w}(X)$, and $N^\circ_k(X)$, $k\ge 2$. {\em (It should be mentioned that the topological center of the semigroup $\beta^\circ(X)$ of free ultrafilters  is empty \cite{Pr98}).}
\end{problem}

\section{Left cancelable elements of $G(X)$}

An element $a$ of a groupoid $S$ is called {\em left cancelable} (resp. {\em right cancelable}) if for any points $x,y\in S$ the equation $ax=ay$ (resp. $xa=ya$) implies $x=y$. In this section we characterize left cancelable elements of the groupoid $G(X)$ over a quasigroup $X$.

\begin{theo}\label{leftcancel} Let $X$ be a quasigroup. An inclusion hyperspace $\F\in G(X)$ is left cancelable in the groupoid $G(X)$ if and only if  $\F$ is a principal ultrafilter.
\end{theo}

\begin{proof} Assume that $\F$ is left cancelable in $G(X)$. First we show that $\F$ contains some singleton.
Assuming the converse, take any point $x_0\in X$ and note that  $F*(X\setminus\{x_0\})=X$ for any $F\in\F$. To see that this equality holds, take any point $a\in X$, choose two distinct points $b,c\in F$ and find solutions $x,y\in X$ of the equation $b*x=a$ and $c*y=a$. Since $X$ is right cancellative, $x\ne y$. Consequently, one of the points $x$ or $y$ is distinct from $x_0$. If $x\ne x_0$, then $a=b*x\in F*(X\setminus\{x_0\})$. If $y\ne x_0$, then $a=c*y\in F*(X\setminus\{x_0\})$.
Now for the inclusion hyperspace $\U=\langle X\setminus\{x_0\}\rangle\ne \min G(X)$, we get $\F\circ\U=\min G(X)=\F\circ\min G(X)$, which contradicts the choice of $\F$ as a left cancelable element of $G(X)$.

Thus $\F$ contains some singleton $\{c\}$. We claim that $\F$ coincides with the principal ultrafilter generated by $c$. Assuming the converse, we would conclude that $X\setminus\{c\}\in \F$. Let $\mathcal A=\langle X\setminus\{c\}\rangle^\perp$ be the inclusion hyperspace consisting of subsets that meet $X\setminus\{c\}$. It is clear that $\A\ne\max G(X)$. We claim that $\F\circ\A=\max G(X)=\F\circ\max G(X)$ which will contradict the left cancelability of $\F$. Indeed, given any singleton $\{a\}\in\max G(X)$, consider two cases: if $a\ne c*c$, then we can find a unique $x\in X$ with $c*x=a$. Since $x\ne c$, $\{x\}\in\A$ and hence $\{a\}=c*\{x\}\in\F\circ\A$. If $a=c*c$, then for every $y\in X\setminus \{c\}$ we can find $a_y\in X$ with $y*a_y=a$ and use the left cancelativity of $X$ to conclude that $a_y\ne c$ and hence $\{a_y\}\in\A$. Then $\{a\}=\bigcup_{y\in X\setminus\{c\}}y*\{a_y\}\in\F\circ\A$.

Therefore $\F=\langle c\rangle$ is a principal ultrafilter, which proves the ``only if'' part of the theorem. To prove the ``if'' part, take any principal ultrafilter $\langle x\rangle$ generated by a point $x\in X$. We claim that two inclusion hyperspaces $\F,\U\in G(X)$ are equal provided $\langle x\rangle\circ\F=\langle x\rangle\circ \U$. Indeed, given any set $F\in\F$ observe that $x*F\in\langle x\rangle\circ\F=\langle x\rangle \circ\U$ and hence $x*F=x*U$ for some $U\in\U$. The left cancelativity of $X$ implies that $F=U\in\U$, which yields $\F\subset\U$. By the same argument we can also check that $\U\subset\F$.
\end{proof}

\begin{problem} Given an (infinite) group $X$ describe left cancelable elements of the subsemigroups $\lambda(X)$, $\Fil(X)$, $N_{<\w}(X)$, $N_k(X)$, $k\ge 2$ (and $G^\circ(X)$, $\lambda^\circ(X)$, $\Fil^\circ(X)$, $N^\circ_{<\w}(X)$, $N^\circ_k(X)$,  for $k\ge 2$). 
\end{problem}

\begin{rem} Theorem~\ref{leftcancel} implies that for a countable Abelian group $X$ the set of left cancelable elements in $G(X)$ coincides with $X$. On the other hand, the set of (left) cancelable elements of $\beta(X)$ contains an open dense subset of $\beta^\circ(X)$, see Theorem 8.34 of \cite{HS}. 
\end{rem}

\section{Right cancelable elements of $G(X)$}

As we saw in the preceding section, for any quasigroup $X$ the groupoid $G(X)$ contains only trivial left cancelable elements. For right cancelable elements the situation is much more interesting. First note that the right cancelativity of an inclusion hyperspace $\F\in G(X)$ is equivalent to the injectivity of the map
$\mu_X\circ G\bar R_\F:G(X)\to G(X)$ considered at the begining of Section~\ref{s2}.
We recall that $\mu_X:G^2(X)\to G(X)$ is the multiplication of the monad $\mathbb G=(G,\mu,\eta)$ while $\bar R_\F:\beta X\to G(X)$ is the Stone-\v Cech extension of the right shift $R_\F:X\to G(X)$, $R_\F:x\mapsto x*\F$. The map $\bar R_\F$ certainly is not injective if $R_\F$ is not an embedding, which is equivalent to the discreteness of the indexed set $\{x*\F:x\in X\}$ in $G(X)$. Therefore we have obtained the following necessary condition for the right cancelability.

\begin{stat}\label{rcp1} Let $X$ be a groupoid. If an incluison hyperspace $\F\in G(X)$ is right cancelable in $G(X)$, then the indexed set $\{x\,\F:x\in X\}$ is discrete in $G(X)$ in the sense that each point $x\F$ has a neighborhood $O(x\F)$ containing no other points $y\F$ with $y\in X\setminus\{x\}$.
\end{stat}

Next we give a sufficient condition of the right cancelability.

\begin{stat}\label{rcp2} Let $X$ be a groupoid. An inclusion hyperspace $\F\in G(X)$ is right cancelable in $G(X)$ provided there is a family of sets $\{S_x\}_{x\in X}\subset \F\cap\F^\perp$  such that $xS_x\cap yS_y=\emptyset$ for any distinct $x,y\in X$.
\end{stat}

\begin{proof} Assume that $\A\circ\F=\mathcal B\circ \F$ for two inclusion hyperspaces $\A,\mathcal B\in G(X)$. First we show that $\A\subset\mathcal B$. Take any set $A\in\A$ and observe that the set $\bigcup_{a\in A}aS_a$ belongs to $\A\circ\F=\mathcal B\circ\F$. Consequently, there is a set $B\in\mathcal B$ and a family of sets $\{F_b\}_{b\in B}\subset\F$ such that $$\bigcup_{b\in B}bF_b\subset \bigcup_{a\in A}aS_a.$$ It follows from $S_b\in\F^\perp$ that $F_b\cap S_b$ is not empty for every $b\in B$.

Since the sets $aS_a$ and $bS_b$ are disjoint for different $a,b\in X$, the inclusion $$\bigcup_{b\in B}b(F_b\cap S_b)\subset \bigcup_{b\in B}bF_b\subset \bigcup_{a\in A}aS_a$$implies $B\subset A$ and hence $A\in\mathcal B$.

By analogy we can prove that $\mathcal B\subset\mathcal A$.
\end{proof}

Propositions~\ref{rcp1} and \ref{rcp2} imply the following characterization of right cancelable ultrafilters in $G(X)$ generalizing a known characterization of right cancelable elements of the semigroups $\beta X$ , see \cite[8.11]{HS}.

\begin{corollary} Let $X$ be a countable groupoid. For an ultrafilter $\U$ on $X$ the following conditions are equivalent:
\begin{enumerate}
\item $\U$ is right cancelable in $G(X)$;
\item $\U$ is right cancelable in $\beta X$;
\item the indexed set $\{x\,\U:x\in X\}$ is discrete in $\beta X$;
\item there is an indexed family of sets $\{U_x\}_{x\in X}\subset\U$ such that for any distinct $x,y\in X$ the shifts $x\,U_x$ and $y\,U_y$ are disjoint.
\end{enumerate}
\end{corollary}

This characterization can be used to show that for any countable group $X$ the semigroup $\beta^\circ(X)$ of free ultrafilters contains an open dense subset of right cancelable ultrafilters, see \cite[8.10]{HS}. It turns out that a similar result can be proved for the semigroup $G^\circ(X)$.

\begin{stat} For any countable quasigroup, the groupoid $G^\circ(X)$ contains an open dense subset of right cancelable free inclusion hyperspaces.
\end{stat}

\begin{proof} Let $X=\{x_n:n\in\w\}$ be an injective enumeration of the countable quasigroup $X$. Given a free inclusion hyperspace $\F\in G^\circ(X)$ and a neighborhood $O(\F)$ of $\F$ in $G^\circ(X)$, we should find a non-empty open subset in $O(\F)$. Without loss of generality, the neighborhood $O(\F)$ is of basic form:
$$O(\F)=G^\circ(X)\cap U_0^+\cap\dots\cap U_n^+\cap U_{n+1}^-\cap\dots\cap U_{m-1}^-$$
for some sets $U_1,\dots,U_{m-1}$ of $X$. Those sets are infinite because $\F$ is free. We are going to construct an infinite set $C=\{c_n:n\in\w\}\subset X$ that has infinite intersection with the sets $U_i$, $i<m$, and such that for any distinct $x,y\in X$ the intersection $xC\cap yC$ is finite. The points $c_k$, $k\in\w$, composing the set $C$ will be chosen by induction to satisfy the following conditions:
\begin{itemize}
\item $c_k\in U_j$ where $j=k \mod m$;
\item $c_k$ does not belong to the finite set $$F_k=\{z\in X:\exists i,j\le k\;\exists l<k\;\;(x_iz=x_jc_l)\}.$$
\end{itemize}
It is clear that the so-constructed set $C=\{c_k:k\in\w\}$ has infinite intersection with each set $U_i$, $i<m$. Since $X$ is right cancellative, for any $i<j$ the set $Z_{i,j}=\{z\in X:x_iz=x_jz\}$ is finite. Now the choice of the points $c_k$ for $k>j$ implies that $x_iC\cap x_jC\subset x_i(Z_{i,j}\cup \{c_l:l\le j\})$ is finite.

Now let $\C$ be the free inclusion hyperspace on $X$ generated by the sets $C$ and $U_0,\dots,U_n$. It is clear that $\C\in O(\F)$ and $C\in\C\cap\C^\perp$. Consider the open neighborhood $$O(\C)=O(\F)\cap C^+\cap (C^+)^\perp$$ of $\C$ in $G^\circ(X)$. 

We claim that each inclusion hyperspace $\A\in O(\C)$ is right cancelable in $G(X)$.
This will follow from Proposition~\ref{rcp2} as soon as we construct a family of sets $\{A_i\}_{i\in\w}\in \A\cap\A^\perp$ such that $x_iA_i\cap x_jA_j=\emptyset$ for any numbers $i<j$. The sets $A_i$, $i\in\w$, can be defined by the formula $A_k=C\setminus F_k$ where 
$$F_k=\{c\in C:\exists i<k\mbox{ with }x_kc=x_i C\}$$ is finite by the choice of the set $C$.
\end{proof} 

\begin{problem} Given an (infinite) group $X$ describe right cancelable elements of the subsemigroups $\lambda(X)$, $\Fil(X)$, $N_{<\w}(X)$, $N_k(X)$, $k\ge 2$ \textup{(}$\lambda^\circ(X)$, $\Fil^\circ(X)$, $N^\circ_{<\w}(X)$, $N^\circ_k(X)$,  for $k\ge 2$\textup{)}. 
\end{problem}

\section{The structure of the semigroups $G(H)$ over finite groups $H$}

In Proposition~\ref{ideall} we have seen that the structural
properties of the finite semigroup $\lambda(\IZ_3)$ has
non-trivial implications for the essentially infinite object
$\lambda^\circ(\IZ)$. This observation is a motivation for more
detail study of spaces $G(H)$ over finite Abelian groups $H$.
In this case the group $H$ acts on $G(H)$ by right shifts:
$$s:G(H)\times H\to G(H),\; s:(\A,h)\mapsto \A\circ h.$$
So we can speak about the orbit $\A\circ H=\{\A\circ h:h\in H\}$
of an inclusion hyperspace $\A\in G(H)$ and the orbit space
$G(H)/H=\{\A\circ H:\A\in G(X)\}$. By $\pi:G(H)\to G(H)/H$ we
denote the quotient map which induces a unique semigroup structure
of $G(H)/H$ turning $\pi$ into a semigroup homomorphism.

We shall say that the semigroup $G(H)$ is {\em splittable} if
there is a semigroup homomorphism $s:G(H)/H \to G(H)$ such that
$\pi\circ s$ is the identity homomorphism of $G(H)/H$. Such a
homomorphism $s$ will be called a {\em section} of $\pi$ and the
semigroup $T(H)=s(G(H)/H)$ will be called a {\em $H$-transversal
semigroup} of $G(H)$. It is clear that a $H$-transversal semigroup
$T(H)$ has one-point intersection with each orbit of $G(H)$.

If the semigroup $G(H)$ is splittable, then the structure of
$G(H)$ can be described as follows.

\begin{stat} If the semigroup $G(H)$ is splittable and $T(H)$ is the transversal
semigroup of $G(H)$, then $T(H)$ is isomorphic to $G(H)/H$ and
$G(H)$ is the quotient semigroup of the product $T(H)\times H$
under the homomorphism $h:T(H)\times H\to G(H)$, $h:(\A,h)\mapsto
\A\circ h$.
\end{stat}

It turns out that the semigroup $G(\IZ_n)$ is splittable for $n\le 3$ and not splittable for $n=5$ (the latter follows from the non-splittability of the semigroup $\lambda(\IZ_5)$ established in \cite{BGN}).  So below we describe the structure of the semigroups
$G(\IZ_n)$ and their transversal semigroup $T(\IZ_n)$ for $n\le
3$. 

For a group $X$ we shall identify the elements $x\in X$ with the ultrafilters they generate. Also we shall use the notations $\wedge$ and $\vee$ to denote the lattice operations $\cap$ and $\cup$ on  $G(X)$, respectively.

{\bf The semigroup $G(\IZ_2)$}. For the cyclic group $\IZ_2=\{e,a\}$ the lattice $G(\IZ_2)$ contains 
four inclusion hyperspaces: $e,a,e\wedge a,e\vee a$, and is shown at the picture:

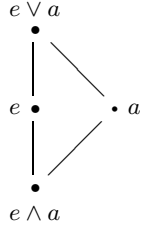
\begin{picture}(100,100)(-150,-20)
\put(0,0){\circle*{3}}
\put(-1,5){\line(0,1){20}}
\put(5,5){\line(1,1){20}}
\put(0,30){\circle*{3}}
\put(-10,28){\small $e$}
\put(35,28){\small $a$}
\put(-1,35){\line(0,1){20}}
\put(0,60){\circle*{3}}
\put(-10,65){\small $e\vee a$}

\put(6,55){\line(1,-1){20}}
\put(30,30){\circle*{2}}
\put(-10,-12){\small $e\wedge a$}
\end{picture}

The semigroup $G(\IZ_2)$ has a unique $\IZ_2$-transversal semigroup $$T(\IZ_2)=\{e\wedge a,e,e\vee a\}$$ with two right zeros: $e\wedge a$, $e\vee a$ and one unit $e$.

{\bf The semigroup $G(\IZ_3)$} over the cyclic group $\IZ_3=\{e,a,a^{-1}\}$ contains 18 elements: 

\centerline{$a\vee e\vee a^{-1}$,}

\centerline{$a\vee a^{-1}$, \ $a\vee e$, \ $e\vee a^{-1}$,}

\centerline{$a\vee (e\wedge a^{-1})$, \ $e\vee (a\wedge a^{-1})$, \ $a^{-1}\vee (a\wedge e)$,}

\centerline{$a,e,a^{-1}$,}

\centerline{$(a\vee  e)\wedge(a\vee a^{-1})\wedge (e\vee a^{-1})$,}

\centerline{$a\wedge (e\vee a^{-1})$, \ $e\wedge (a\vee a^{-1})$, \ $a^{-1}\wedge (a\vee e)$,}

\centerline{$a\wedge a^{-1}$, \ $a\wedge e$, \ $e\wedge a^{-1}$,}

\centerline{$a\wedge e\wedge a^{-1}$}

\noindent divided into 8 orbits with respect to the action of the group $\IZ_3$. 

The semigroup $G(\IZ_3)$ has 9 different $\IZ_3$-transversal semigroups one of which is drawn at the picture:

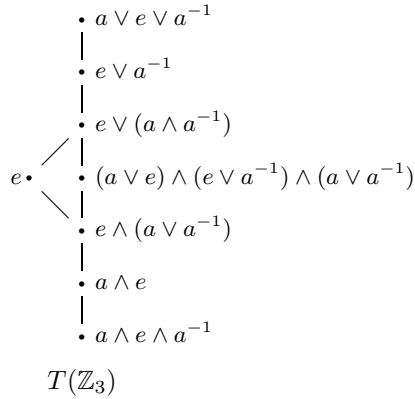
\begin{picture}(100,170)(-30,-30)
\put(117,-20){$T(\IZ_3)$}

\put(130,0){\circle*{2}}
\put(130,5){\line(0,1){10}}
\put(135,-2){\small $a\wedge e\wedge a^{-1}$}
\put(130,20){\circle*{2}}
\put(130,25){\line(0,1){10}}
\put(135,18){\small $a\wedge e$}
\put(130,40){\circle*{2}}
\put(130,45){\line(0,1){10}}
\put(125,45){\line(-1,1){10}}
\put(135,38){\small $e\wedge (a\vee a^{-1})$}
\put(110,60){\circle*{2}}
\put(130,65){\line(0,1){10}}
\put(103,58){\small $e$}
\put(130,60){\circle*{2}}
\put(130,65){\line(0,1){10}}
\put(135,58){\small $(a\vee e)\wedge (e\vee a^{-1})\wedge(a\vee a^{-1})$}
\put(130,80){\circle*{2}}
\put(130,85){\line(0,1){10}}
\put(135,78){\small $e\vee (a\wedge a^{-1})$}
\put(125,75){\line(-1,-1){10}}
\put(130,100){\circle*{2}}
\put(130,105){\line(0,1){10}}
\put(135,98){\small $e\vee a^{-1}$}
\put(130,120){\circle*{2}}
\put(135,118){\small $a\vee e\vee a^{-1}$}

\end{picture}

The semigroup $G(\IZ_3)$ has 3 shift-invariant inclusion hyperspaces which are right zeros: $a\wedge e\wedge a^{-1}$, $a\vee e \vee a^{-1}$ and $(a\vee e)\wedge (e\vee a^{-1})\wedge(a\vee a^{-1})$. Besides right zeros $G(\IZ_3)$ has 3 idempotents:
$e$, $e\vee (a\wedge a^{-1})$ and $e\wedge (a\vee a^{-1})$. The element $e$ is the unit of the semigroup $G(\IZ_3)$.

The complete information on the structure of the $\IZ_3$-transversal semigroup $T(\IZ_3)$ (which is isomorphic to the quotient semigroup $G(\IZ_3)/\IZ_3$) can be derived from the Cayley table 

\hskip80pt\begin{tabular}{c|ccccccc}
$\circ$ &$x_{-3}$&$x_{-2}$&$x_{-1}$&$x_0$& $x_1$& $x_2$&$x_3$\\
\hline
$x_{-3}$&$x_{-3}$&$x_{-3}$&$x_{-3}$&$x_0$& $x_0$& $x_0$&$x_3$\\

$x_{-2}$&$x_{-3}$&$x_{-3}$&$x_{-2}$&$x_0$& $x_0$& $x_1$&$x_3$\\

$x_{-1}$&$x_{-3}$&$x_{-3}$&$x_{-1}$&$x_0$& $x_0$& $x_2$&$x_3$\\

$x_0$   &$x_{-3}$&$x_{-3}$& $x_0$  &$x_0$& $x_0$& $x_3$&$x_3$\\

$x_1$   &$x_{-3}$&$x_{-2}$& $x_0$  &$x_0$& $x_1$& $x_3$&$x_3$\\

$x_2$   &$x_{-3}$&$x_{-1}$& $x_0$  &$x_0$& $x_2$& $x_3$&$x_3$\\

$x_3$   &$x_{-3}$&$x_0$   & $x_0$  &$x_0$& $x_3$& $x_3$&$x_3$
\end{tabular}
\smallskip

\noindent of its  linearly ordered subsemigroup $T(\IZ_3)\setminus\{e\}$ having  with 7-elements:
$$
\begin{aligned}
x_{-3}=&\;e\wedge a\wedge a^{-1},\\
x_{-2}=&\;e\wedge a,\\
x_{-1}=&\; e\wedge (a\vee a^{-1}),\\
x_0=&\; (e\vee a)\wedge (e\vee a^{-1})\wedge(a\vee a^{-1}),\\
x_1=&\;e\vee (a\wedge a^{-1}),\\
x_2=&\;e\vee a,\\
x_3=&\;e\vee a\vee a^{-1}.
\end{aligned}
$$

\section{Acknowledgments}

The author express his sincere thanks to Taras Banakh and Oleg Nykyforchyn for help during preparation of the paper and also to the referee for inspiring criticism.

\end{document}